# A System-Dynamic Based Simulation and Bayesian Optimization for Inventory Management

Sarit Maitra
Sarit.maitra@gmail.com; sarit.maitra@alliance.edu.in
Alliance University, Bangalore, India

*Abstract*— Inventory management is a fundamental challenge in supply chain management. The challenge is compounded when the associated products have unpredictable demands. This study proposes an innovative optimization approach combining system-dynamic Monte-Carlo simulation and Bayesian optimization. The proposed algorithm is tested with a real-life, unpredictable demand dataset to find the optimal stock to meet the business objective. The findings show a considerable improvement in inventory policy. This information is helpful for supply chain analytics decision-making, which increases productivity and profitability. This study further adds sensitivity analysis, considering the variation in demand and expected output in profit percentage. This paper makes a substantial contribution by presenting a simple yet robust approach to addressing the fundamental difficulty of inventory management in a dynamic business environment.

## I. Introduction

Inventory management (IM) stands as a cornerstone in supply chain management and contributes significantly to a company's financial performance. It is a complex task, and the complexities are multiplied when combined with unpredictability. The global supply chain network faces immense challenges because of the competitive business environment, changing customer demand patterns, a wide variety of product types, shortened product lifecycles, and dynamic global marketplaces. Improper management can lead to major disruptions, e.g., stockouts, overstocking, backorders, etc. Moreover, a modern supply chain often employs multiple tiers of stocking locations to reduce logistics costs, inventory costs, and sales opportunity losses. This gives rise to additional operational risk-related difficulties.

Operational risks involve daily disturbances in inventory operations like lead-time and demand fluctuations. The disruption risks often involve low-frequency but high-impact events ([1], [2]). The recent pandemic is an example of such disruption which led to material shortages, delivery delays, and performance degradation, affecting revenue, service levels, and productivity ([3], [4], [5]). Keeping in view of such events, recent researchers have turned their focus to uncertain demand modeling, which provides a flexible and principled method for probabilistic modelling and inference (e.g., [6], [7], [8], [9], [10], [11], [12]). The optimal order quantity and reorder point in todays' dynamic business scenario is a critical aspect for business success. Here optimization process plays an indispensable role for applying appropriate inventory policy and resilient supply chain. The main goal is profit maximization through better customer service, minimize overall costs by reducing stockouts and revenue losses and maintain service level agreements [13].

The use of simulations and advanced modeling techniques are critical for a dynamic representation of real-life challenges, enabling decision-makers to test different scenarios and strategies in a controlled environment (e.g., [14], [15], [16], [17], [18], [19] etc.). Simulations involve computer based model, is a valuable tool for establishing effective policies or testing the results of mathematical models in the context of uncertain demand [20]. Numerous studies have highlighted the effectiveness of simulations (e.g., [21], [22], [23], [24], [25], [26], [27] [28], [29], [30], [31], [32] etc.). The dynamic nature of multi-echelon inventory management, along with factors such as lead times, variability in demand, and inventory holding costs, makes it challenging to derive precise analytical solutions. These often require to be evaluated using a simulator. Though simulations are useful to replicate real-world scenarios and study probabilistic demand patterns, they are most effective when combined with optimization process.

On top of the simulation model, an optimization process is necessary to determine the best variables given a predefined objective function. Genetic algorithms are commonly employed to identify ideal IM configurations. While these methods are effective, they are complex and computationally expensive due to the need to run the simulator for multiple scenarios and evaluate the objective function under numerous parameter settings. We present a hybrid approach that combines Monte Carlo simulations (MCS) and Bayesian optimization (BayesOpt) for determining the ideal IM configuration and evaluating its efficiency using computational trials.

While the existing literature has made substantial strides in addressing the challenges of stochastic demand and optimizing inventory policies, it is imperative to acknowledge that the field continues to evolve rapidly. This work aims to improve IM decisions in uncertain environments. It uses an analytical and simulation-based framework for a deeper understanding of inventory policies and optimization techniques. To be precise, an inventory model that specifically takes into account uncertainty over both stationary and non-stationary demand models is proposed in this study. This comprises constraint-based simulation modeling, which allows for experimentation with different scenarios and policies. It is further integrated with probabilistic BayesOpt to optimize the order quantity to maximize profit. The primary contribution of this work is that it provides system dynamics-based optimization modeling for inventory management for retail supplies. The proposed model can be viewed as a helpful mechanism to enhance service levels, safety stocks, order quantity, and the associated overall profits related to inventory management in a retail supply chain.

## II. METHODOLOGY

The objective here is to identify the optimal inventory level that maximizes profit while ensuring high service levels (95%). We employed a hybrid approach where, in the first step, by applying Monte Carlo simulations (MCS), we analyze different demand scenarios and make informed decisions to optimize inventory levels. In the second step, we integrate BayesOpt to further optimize the output from MCS.

MCS is employed to estimate the expected inventory cost under uncertainty. The inventory-optimization problem for a given policy:

$$\min \frac{1}{n}\sum_1^n Inv(\alpha, s_i), s_i \rightarrow scenario\ i\ \epsilon\ \{1, 2, \ldots, n\}$$

$$\text{Subject to,}\ \min_{i,j} Fr_j(\alpha, s_i) \geq \beta$$

where, $Inv(\alpha, s_i) = Inventory\ cost$, $Fr_j(\alpha, s_i)$ = $cumulative\ demand\ fulfill\ rate\ at\ market\ site\ j$

$\alpha = inventory\ policy\ parameter$, $\beta = minimum\ demand\ fulfill\ rate\ required\ at\ each\ j$

The constraint here is to ensure that the demand-fulfill rate at any market site is not lower than $\beta$ for any scenario.

Two proposed methods for handling this constraint are penalty-based BayesOpt and constraint BayesOpt [33]. The former introduces a high penalty cost when a parameter setting falls into an infeasible region, while the latter approximates the constraint function using a Gaussian Process Regression (GPR) surrogate model are considered for this work.

For a function $f$ depending on $\theta$, a vector of unknown parameters is estimated based on data $x$,

$$E(f(\theta)|x) \approx \frac{1}{n}\sum_{s=1}^{n} f(\theta^s) \qquad (1)$$

$\theta^s$ represents a simulated sample from the posterior probability distribution of $\theta$ and $n$ is the total number of simulated samples. We need to determine the expected value of function $(f(\theta))$ given the data $(x)$. This is determined by taking the average of $f(\theta^s)$ over $n$ simulated samples and involves GPR with conditioning function to approximate the objective function and EI is used as acquisition function.

$$f(x) \sim GP(m(x), k(x_i, x_j)) \qquad (2)$$

$m(x)$ and $k(x_i, x_j)$) are the mean and covariance functions (kernel function) and denoted as:

$$m(x) = E(f(x))$$
$$k(x_i, x_j) = E\left\{[f(x) - m(x)][f(x_j) - m(x_j)]^T\right\}$$
$$k_{Matern}(x_i, x_j) = \frac{1}{\Gamma(v)2^{v-1}}\left(\frac{\sqrt{2v}}{l}d(x_i, x_j)\right)^v K_v\left(\frac{\sqrt{2v}}{l}d(x_i, x_j)\right)$$
$$(3)$$

Matérn kernel [34] introduces flexibility in modeling, capture pattern, and smoothness in the underlying function. The conditioning function enhances the adaptability of the optimization algorithm. As more data As it accumulates more data, it modifies the mean function $m(x)$ based on available data:

$$m_c(x) = m(x) + \mathbb{C}(x) \qquad (4)$$

here $m_c(x)$ is the conditioned mean function. The conditioning function defined as:

$$\mathbb{C}(x) = \alpha * u(x) \qquad (5)$$

where $\alpha$ = scaling factor and $u(x)$ captures the influence of observed data on the mean function. The combined objective function is:

$$f_{cmb}(x) = m_c(x) + f(x) \qquad (6)$$

The goal is to find the optimal input $x^*$ that maximizes or minimizes the combined objective function.

$$x^* = arg\ max_x f_{cmb}(x)) \qquad (7)$$

This means we are looking for the input that gives us the best balance between our prior expectations (captured by $m(x)$) and the observed data (captured by $\mathbb{C}(x)$), while optimizing the original function $f(x)$. Eq. 8 presents the EI function, which calculates the expected improvement in the objective function value if a new evaluation is performed at input x.

$$EI(x) = \int_{-\infty}^{f_{min}} (f_{min} - f_{cmb}(x))^+ p(f_{cmb}(x)|D) df_{cmb}(x)$$

$$(8)$$

Here, $f_{min}$ = minimization function, $p(f_{cmb}(x)|D)$ = posterior probability distribution of $f_{cmb}(x)$ given the observed data D, and $(\cdot)^+$ = the positive part of the function. In practical terms, this approach seeks the input that optimally balances prior expectations and observed data while optimizing the original objective function $f(x)$.

The system dynamics are embedded in the formulation through the conditioning function's adjustment of the mean function (Eq. 5) and the subsequent integration of the conditioned mean function into the combined objective function (Eq, 8). This allows for the representation of evolving relationships and feedback loops within the system, capturing its dynamic behavior over time.

## III. DATA DESCRIPTION

For empirical analysis, we collected 365 days of sales data to model the demand for four distinct products. Table 1 provides a summary of each product based on past sales data (365 days). This analysis is about uncertainties caused by unpredictable purchase patterns identified by average daily demand and deviations in order. All the products have the same holding costs, which simplifies the cost structure. Pr B is high in demand and gets sold every day (p = 1), and the mean order size is 649 (648.55). The PrD is being purchased 23% of the time, with a mean order quantity of around 150 (150.06). From a system configuration perspective, we used an Intel(R) Core (TM) i5-8265U CPU at 1.60 GHz, Google Cloud, and Python v3.10.12.

TABLE I. STATISTICAL SUMMARY OF THE DEMAND DATA

| Description | Pr A | Pr B | Pr C | Pr D |
|---|---|---|---|---|
| Purchase Cost (PC) | 12 | 7 | 6 | 37 |
| Lead Time (LT) | 9 | 6 | 15 | 12 |
| Order Quantity | 0.57 | 0.05 | 0.53 | 1.05 |
| Selling Price (SP) | 16.10 | 8.60 | 10.20 | 68.00 |
| Starting Stock | 2,750 | 22,500 | 5,200 | 1,400 |
| Mean Order ($\mu$) | 78.33 | 648.55 | 141.61 | 35.67 |
| Std Dev (σ) | 55.08 | 26.48 | 95.96 | 63.98 |
| Order Cost (OC) | 1,000 | 1,200 | 1,000 | 1,200 |
| Holding Cost (HC) | 20 | 20 | 20 | 20 |
| Probability of demand on any day (P) | 76% | 100% | 70% | 23% |
| Demand Lead | 705 | 3891 | 2266 | 785 |
| Total demand | 28,670 | 237,370 | 51,831 | 13,056 |
| Maximum demand | 214 | 718 | 267 | 156 |
| Safety Stock (SS) | 185 | 107 | 199 | 18 |
| Reorder Ooint (ROP) | 1116 | 3998 | 3224 | 1819 |
| Economic Order Quantity (EOQ) | 1693 | 5337 | 2277 | 1252 |

The starting stock is calculated by the average demand during the lead time plus the safety stock. Safety stock is calculated by multiplying the standard deviation of lead time demand by a safety stock coefficient. PrA has a lead time of 9 days, and an average of 705 orders are expected over that time, and likewise for the other products. This needs to be considered while placing the initial order to avoid order loss.

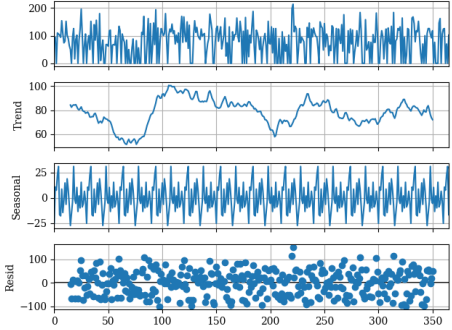

Fig 1. Seasonal-Trend decomposition using LOESS (STL) for Pr A

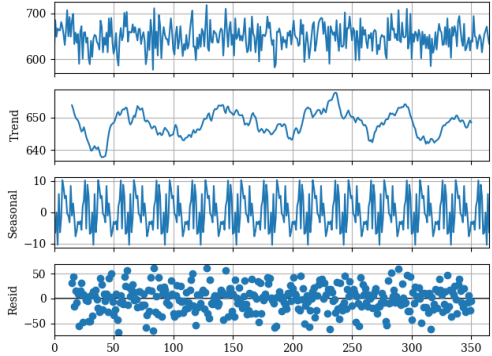

Fig 2. Seasonal-Trend decomposition using LOESS (STL) for Pr B

STL decomposition is applied to check for any recurring temporal pattern that exists in the data. Figs. 1, 2, 3, and 4 display the plots with a frequency of 30 days. We can observe seasonal oscillation from these plots. The remainder looks noisy and is lacking a particular pattern. Residual plots indicate high variability, which indicates unpredictable demands.

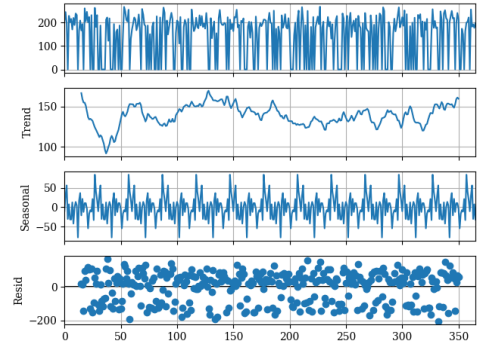

Fig 3. Seasonal-Trend decomposition using LOESS (STL) for Pr C

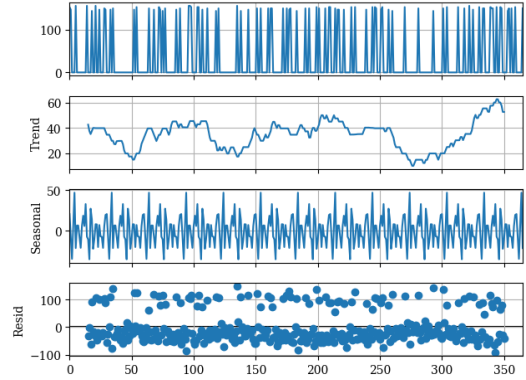

Fig 4. Seasonal-Trend decomposition using LOESS (STL) for Pr D

From Table I, we see uncertainties related to demand variability, lead time variability, and other stochastic factors exist in the historical data. The calculations of SS (safety stock) and ROP (reorder point) are shown below. The lead time (LT) refers to the time interval between placing an order and receiving it. The demand during this LT is uncertain and can be modeled with a probability distribution. Eq. 9 helps in quantifying the variability in demand during LT.

$$\sigma_{demand_{leadtime}} = \sigma_{demand}[i] * \sqrt{LT[i]} \quad (9)$$

$$SS[i] = z_{score} * \sigma \quad (10)$$

$$\mu_{demand_{leadtime}} = \mu[i] * LT[i] \quad (11)$$

$$ROP[i] = \mu_{demand_{leadtime}} + SS[i] \quad (12)$$

$$EOQ = \sqrt{\frac{2 * Annual\ Demand * OC}{HC\ /\ unit}} \quad (13)$$

$$Total\ Annual\ Cost = [(Annual\ Demand\ /\ EOQ) * OC] + [(EOQ\ /\ 2) * HC] \quad (14)$$

However, the EOQ makes a number of assumptions that are not true for every business. It assumes that the unit price, purchase costs, and rate of demand are all constant. Thus, we employ Monte Carlo simulations to model the variability and uncertainty that are not captured in the deterministic EOQ model.

IV. MONTE-CARLO SIMULATIONS

Simulations are used to calculate the appropriate inventory policy for the given dataset.

$$Revenue = SP_i * Units\ sold \quad (15)$$

$$Annual\ profit = Revenue - (PC + OC + HC) \quad (16)$$

$$= SP_i \sum_{t=1}^{365} S_{i,t} - \left[\left(\frac{20v_i}{365}\right)\sum_{t=1}^{365} I_{i,t} + N_i C_{o,i} + \sum_{t=1}^{365} c_i p_{i,t}\right] \quad (17)$$

We perform MCS over a 365-day period to track inventory levels, orders, and profits. The profit is calculated considering revenue, costs, and inventory. The MCS involves running multiple scenarios to account for the variability and uncertainty in demand. Here, we assume that the daily demand follows a normal distribution with parameters determined by the mean and standard deviation of the demand data for each product. This way, modeling demand provides a more realistic representation of demand variability.

A range is specified (starting point as ROP) to conduct a comprehensive performance analysis of each product under different order quantity scenarios. By evaluating multiple order quantities within the range, the model identifies the required order quantity that maximizes the expected profit while considering factors such as profit variability and the risk of lost orders. We apply the same range for both of these policies, incrementing by 10. Table II displays the output.

TABLE II. OUTPUT OF SIMULATION.

| Products | Order quantity | Mean profit | Profit std dev | Mean lost orders |
|---|---|---|---|---|
| Pr A | 2,066 | 88,640 | 3,511 | 6% |
| Pr B | 5,998 | 274,822 | 212 | 46% |
| Pr C | 4,124 | 167,025 | 6,606 | 5% |
| Pr D | 1,829 | 315,245 | 36,516 | 4% |

*Time taken: 193.35 seconds

The standard deviation of profits corresponding to the optimal order point provides a measure of the variability of the expected profit. The proportion of lost orders is relative to the total demand and provides insights into how well the model handles demand fluctuations.

## V. BAYSEAN OPTIMIZATION

We further conduct a comprehensive analysis by employing BayesOpt to optimize the profit based on the optimized quantity derived from MCS. BayesOpt adaptively selects new parameter settings to evaluate based on the current surrogate model and acquisition function. By leveraging information from previous simulations, the process explores the parameter space, focusing on areas with high uncertainty or potential for improvement.

TABLE III. OPTIMIZED OUTPUT.

| Products | Order quantity | Mean profit | Profit std dev | Mean lost orders |
|---|---|---|---|---|
| Pr A | 2,197 | 417,804 | 16,951 | 14% |
| Pr B | 7,358 | 1,237,303 | 1,323 | 38% |
| Pr C | 2,688 | 466,946 | 17,895 | 36% |
| Pr D | 1,585 | 818,974 | 83,393 | 21% |

* Time taken: 158.79 seconds

Compared with MCS, Table III displays how each factor is changed following BayesOpt. The mean profit for all goods appears to have increased significantly, despite the differences in profit standard deviation and mean lost orders.

Fig. 5 displays the inventory distribution for each product based on the optimization process. Each line in the plot represents a different simulation run, showing the changes in the inventory over time for each product under different scenarios. A high level of fluctuation indicates higher variability in demand or lead time. Product 4 (Pr D) has higher demand variability compared to Product 2 (Pr B). PrB has stable and consistent demand. Denser lines indicate more fluctuations and variability in the inventory levels during the simulation runs. This suggests that the order quantity optimization for Pr. D is not effective in maintaining stable inventory levels, leading to more frequent changes in stock levels over time. This is also evident from the high standard. deviation of Pr D (83393). This is in line with the original data.

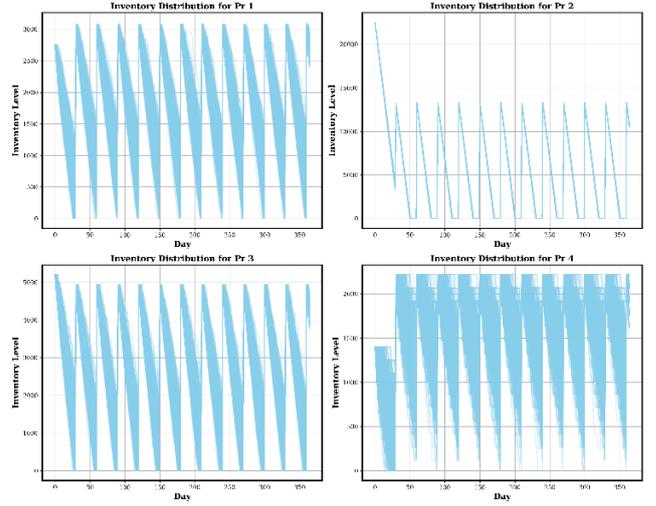

Fig. 5. Inventory distribution for each Product

## VI. EMPIRICAL ANALYSYS

The analysis shows that the best order quantities (annual) for all four products are: [2197, 7358, 2688, 1585] and the average expected profit is: [417804, 1237303, 466946, 818974]. The experimental results indicate that the hybridized Bayes Opt efficiently optimized inventory policy in terms of both optimality and computational efficiency. The simulation, lasting 365 days, covers demand and inventory dynamics for a year, considering daily demand, ordering decisions, and inventory management. Iterations are performed 1000 times for each product to estimate metrics like mean profit, profit standard deviation, and mean lost orders percentage.

### A. Sensitivity Analysis

For each product, we vary each parameter by $+/-$ 10% and 20% to validate the corresponding changes in subsequent parameters. Table V displays the result. We do find that the changes in order quantity do not have a substantial impact on the subsequent parameters such as mean profit, profit standard deviation, and mean lost orders.

TABLE IV. SENSITIVITY ANALYSIS BY VARYING EACH VARIABLE ONE AT A TIME.

| Parameters | Variation in Order Qty | PrA | PrB | PrC | PrD |
|---|---|---|---|---|---|
| Mean profit | +10%/+20% | 417,804 | 1,237,303 | 466,946 | 828,974 |
| | −10%/−20% | | | | |
| Profit Std Dev | +10%/+20% | 16951 | 1323 | 17895 | 83393 |
| | −10%/−20% | | | | |
| Mean Lost Orders | +10%/+20% | 14% | 38% | 36% | 21% |
| | −10%/−20% | | | | |

Fig. 2 displays how changes in the mean profit values for each product affect the overall profit. The height of each bar indicates the corresponding mean profit value for the scenario. We see that Pr B has a significant impact on overall profit, and this product is critical where small adjustments in mean profit can lead to notable changes in overall profit. This is followed by Pr D. Here, we have considered the profit variation range with 5 values evenly spaced between $+/-$ 50,000 range of the mean profit.

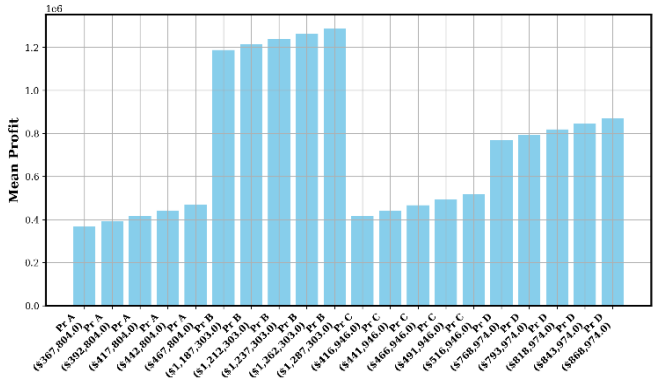

Fig. 6. What-if analysis of profit sensitivity.

Here we observe that the hybridization approach offered a reliable and efficient framework. This aided in managing uncertainty and stochasticity in the optimization process, guided the search for optimal solution. BayesOpt is used to simulate and optimize the order quantity for each product using the same objective function.

Finally, the same method is subsequently applied, utilizing all the factors in the original dataset in Table I (such as purchase cost, lead time, order quantity, selling price, starting stock, mean order, standard deviation, order cost, holding cost, probability of demand on any day, and demand lead) to ensure robust optimization. Table V displays the optimized output considering all the parameters.

TABLE V. OPTIMIZED OUTPUT CONSIDERING ALL PARAMETERS.

| Products | Order quantity | Mean profit | Profit std dev | Mean lost orders |
|---|---|---|---|---|
| Pr A | 2,264 | 416,933 | 17,207 | 7% |
| Pr B | 9,999 | 1,484,243 | 1,327 | 26% |
| Pr C | 2,782 | 467,910 | 17,891 | 24% |
| Pr D | 1,459 | 823,116 | 81,231 | 18% |

*Time taken: 208.83 seconds

Some differences in order quantities and mean lost orders can be seen (compared to Table III), but the differences in mean profit and profit standard deviation are relatively small. We performed a one-way ANOVA test to determine if there were significant differences between the two sets (Tables III and V).

TABLE VI. ANOVA TEST OUTPUT.

| Products | f-stats | p-value |
|---|---|---|
| Order Quantity | 0.079 | 0.787 |
| Mean Profit | 0.040 | 0.846 |
| Profit Std Dev | 0.0003 | 0.985 |
| Mean Loas Orders | 1.386 | 0.283 |

We see from Table VI that the p-values are higher than the significance level of 0.05. Therefore, we fail to reject the null hypothesis, indicating that there is no significant difference between the two sets of products in terms of order quantity, mean profit, profit standard deviation, or mean lost orders. Moreover, constraints helped to develop workable and realistic solutions, which improved the relevance and usefulness of the optimization outcomes in real-world scenarios.

### B. Limitations & Future Directions

The performance of BayesOpt is influenced by hyperparameters like kernel function, acquisition function, and prior distributions. Suboptimal hyperparameter settings may lead to poor convergence behavior. Moreover, BayesOpt assumes that the objective function is a black box with unknown analytical form and derivatives. While this makes it applicable to a wide range of problems, it also limits the ability to exploit problem-specific structures or domain knowledge for optimization. Local optima can cause the optimization to stall, particularly when dealing with non-convex or multimodal objective functions. Global optimization techniques in general may solve this issue; however, they are computationally expensive to operate.

In the future, we shall aim to go for a dynamic, adaptable, and scalable BayesOpt algorithm that can handle high-dimensional input spaces in large-scale optimization situations. This may include different computing environments, e.g., distributed computing, parallelization, or optimization-accelerating approximation approaches.

## VII. CONCLUSION

This paper addresses the inherent uncertainty in stationary or non-stationary demand models by presenting a hybrid optimization technique that combines Bayesian optimization and Monte Carlo simulations. The proposed approach can be used to solve a global optimization problem with multimodal, noisy, and complex objective functions. Even in difficult optimization landscapes, the proposed technique effectively finds optimal solutions by iteratively exploring the parameter space and improving the surrogate model. This high degree of flexibility and adaptability allows researchers and practitioners to customize this simulation model, surrogate model, and optimization technique to match the unique requirements of their situation. This hybrid approach estimates the relationship between the input parameters and the objective function's value. Through this process of sequentially selecting, evaluating, and updating, Bayesian optimization converges to the optimal set of parameter values with fewer evaluations compared to brute-force or random search methods.